\documentclass[11pt]{article}
\usepackage{amssymb,amsmath,latexsym}
\newtheorem{theorem}{Theorem}
\newtheorem{definition}{Definition}
\newtheorem{lemma}{Lemma}
\newtheorem{remark}{Remark}
\newenvironment{proof}{\noindent{\bf Proof: }}{$\hfill \Box$ \vspace{10pt}}

\title{Exploring the zeros of real self-reciprocal polynomials by Chebyshev polynomials\footnote{This work is supported by grant \#2016/02700-8, S\~ao Paulo Research Foundation (FAPESP).}}

\author{Vanessa Botta\\
Universidade Estadual Paulista (Unesp)\\
 Faculdade de Ci\^encias e Tecnologia\\
  C\^ampus de Presidente Prudente - SP, Brazil \\
E-mail: botta@fct.unesp.br}

\date{\today}

\begin{document}

\maketitle

\begin{abstract}In this paper we present some classes of real self-reciprocal polynomials with at most two zeros outside the unit circle which are connected with a Chebyshev quasi-orthogonal polynomials of order one. We investigated the distribution, simplicity and monotonicity of their zeros around the unit circle and real line.\\

\end{abstract}




\section{Introduction}

Let the polynomial $P(z)=\displaystyle{\sum_{i=0}^{n}}p_iz^i$, $p_i\in \mathbb{C}$, $p_n\neq 0$, such that
$P(z)=\displaystyle{z^n P\left(1/z\right)}$. Then $P(z)$ is said to be self-reciprocal polynomial (\cite{Joyner,KimPark}). Self-reciprocal polynomials can be found in the literature under various names such as reciprocal (\cite{Andrews,LAK,Losonczi}) and palindromic (\cite{Blake,Fulman}). It is clear that if $P(z)$ is a self-reciprocal polynomial then $p_i=p_{n-i}$, for $i=0,1,...,n$.  

In this paper we consider the real self-reciprocal polynomials (i.e., $p_i\in \mathbb{R}$) and we present some classes of real self-reciprocal polynomials with at most two zeros located outside the unit circle. We show that these classes of polynomials are  connected with the linear combination of Chebyshev polynomials and we focus in the case that the linear combination of Chebyshev polynomials is a Chebyshev quasi-orthogonal polynomial of order one.

\section{Preliminary results}
\label{section_preliminaries}

The theory of Chebyshev polynomials  is a classical topic and was many explored in the literature, see \cite{Mason,Rivlin}. From \cite{Mason} we have:

\begin{definition}\label{secondkind}
The Chebyshev polynomials $T_n(x)$, $U_n(x)$, $V_n(x)$ and $W_n(x)$ of first, second, third and fourth kinds, respectively, are polynomials in $x$ of degree $n$ defined  by  
\begin{eqnarray*}
&&T_{n}(x)= \cos{n\theta}, \mbox{  } U_{n}(x)= \frac{\sin{(n+1)\theta}}{\sin \theta}, \mbox{ } V_{n}(x)= \frac{\cos{\left(n+\frac{1}{2}\right)\theta}}{\cos \frac{1}{2}\theta} \mbox{ and }\\
&&W_{n}(x)= \frac{\sin{\left(n+\frac{1}{2}\right)\theta}}{\sin \frac{1}{2}\theta},
\end{eqnarray*}
when  $x=\cos\theta$.
\end{definition}


The sequence of polynomials $\displaystyle{\{T_n(x)\}_{n=0}^{\infty}}$, $\displaystyle{\{U_n(x)\}_{n=0}^{\infty}}$, $\displaystyle{\{V_n(x)\}_{n=0}^{\infty}}$ and $\displaystyle{\{W_n(x)\}_{n=0}^{\infty}}$ are orthogonal in $[-1,1]$ with respect to weight functions \linebreak $w(x)=(1-x^2)^{-1/2}$, $w(x)=(1-x^2)^{1/2}$, $w(x)=(1+x)^{1/2}(1-x)^{-1/2}$ and $w(x)=(1+x)^{-1/2}(1-x)^{1/2}$, respectively, and satisfy the same three term recurrence relation with different initial conditions. For example, $\displaystyle{\{T_n(x)\}_{n=0}^{\infty}}$ satisfy
\begin{eqnarray}\label{TTRR}
T_n(x)=2xT_{n-1}(x)-T_{n-2}(x), \mbox{  } n=2,3,\ldots,
\end{eqnarray}
with $T_0(x)=1$ and $T_1(x)=x$.


In this section we develop all the necessary preliminaries to show the main result about the  distribution, simplicity and monotonicity of the zeros of some classes of real self-reciprocal polynomials, presented in Section 3.


From \cite{Mason}, the following relations hold:
\begin{eqnarray}
2T_n(x)=U_n(x)-U_{n-2}(x) \label{relationU}, \\
2T_n(x)=V_n(x)+V_{n-1}(x) \label{relationV} \mbox{ and }\\
2T_n(x)=W_n(x)-W_{n-1}(x) \label{relationW}. 
\end{eqnarray}

Let $\displaystyle{P(z)=\sum_{j=0}^{2n}p_{j}z^{j}}$, $p_{2n}\neq 0$, be a self-reciprocal polynomial of degree $2n$, where $p_j \in \mathbb{R}$, $j=1,\ldots,2n$. Hence
\begin{eqnarray} \label{sumtrans}P(z)=\sum_{j=0}^{2n}p_{j}z^{j}=2z^{n}\left[\frac{p_{2n}}{2}\left(z^{n}+ \frac{1}{z^{n}}\right)+\ldots+\frac{p_{n+1}}{2}\left(z+ \frac{1}{z}\right)+\frac{p_{n}}{2}\right].
\end{eqnarray}


Considering the transformation
\begin{equation}\label{changevariable}
\displaystyle{\frac{1}{2}\left(z+\frac{1}{z}\right)}=x,
\end{equation}
from relation \eqref{TTRR} and basic manipulations
follows that  $\displaystyle{\frac{1}{2}\left(z^{j}+\frac{1}{z^{j}}\right)}=T_{j}(x)$, $j=0,1,2,\ldots$.

From (\ref{sumtrans}),
\begin{eqnarray} \label{sumtrans2}
P(z)&=&2z^{n}\left[p_{2n}T_n(x)+ p_{2n-1}T_{n-1}(x)+\ldots+p_{n+1}T_1(x)+\frac{p_n}{2}T_0(x)\right]\nonumber \\
&=&2z^nC_n(x).
\end{eqnarray}

Observe that the behaviour of the zeros of real self-reciprocal polynomials is connected with the behaviour of the zeros of a linear combination of Chebyshev polynomials of first kind. In fact, the linear combinations of orthogonal polynomials represents a quasi-orthogonal polynomial, which is a theory very explored by \cite{Chihara},  where we find the following result:
\begin{theorem}[\cite{Chihara}]
Let $\{Q_n(x)\}_{n=0}^{\infty}$ be the family of orthogonal polynomials on $[a, b]$ with respect to a positive weight function $w(x)$. A necessary and sufficient condition for a polynomial $R_n(x)$ of degree $n$ to be quasi-orthogonal of order $r$ on $[a, b]$ with respect to $w(x)$ is that
$$R_n(x) = c_0Q_n(x) +c_1Q_{n-1}(x)+ \cdots +c_rQ_{n-r}(x),$$
where the $c_i$'s are numbers which can depend on $n$ and $c_0c_r\neq 0$.
\end{theorem}

When $r\geq 1$, the following result provides the number of zeros  of the quasi-orthogonal polynomial $R_n(x)$ in $(a, b)$.

\begin{theorem}[\cite{Shohat}]  \label{zeros}
If $R_n(x)$ is quasi-orthogonal of order $r$ on $[a, b]$ with respect to a positive weight function $w(x)$, then at least $n-r$ distinct zeros of $R_n(x)$ lie in the interval $(a, b)$.
\end{theorem}

In the literature there are some results on the behaviour of the zeros of Chebyshev quasi-orthogonal polynomials (see \cite{Alfaro, Brezinski}). In this manuscript we analyse the classes of real self-reciprocal polynomials that are connected with  Chebyshev quasi-orthogonal polynomials of order one. Hence, we use a result about  the properties of the zeros of
\begin{eqnarray*}
R_n(x)=Q_n(x)-cQ_{n-1}(x), c\neq 0,
\end{eqnarray*}
 presented in \cite{BracDimRanga} and \cite{Brezinski}.

Let $\{Q_n(x)\}_{n=0}^{\infty}$ be the family of orthogonal polynomials on $[a, b]$ with respect to a positive weight function $w(x)$. Let $a<x_{n-1,1}<x_{n-1,2}<\ldots<x_{n-1,n-1}<b$ be the zeros of $Q_{n-1}$ and $a<x_{n,1}<x_{n,2}<\ldots<x_{n,n}<b$ those of $Q_n$. We consider
$\displaystyle{f_n(x)=\frac{Q_n(x)}{Q_{n-1}(x)}.}$

Suppose, without loss of generality, that the leading coefficients of all polynomials $Q_n$ have the same sign. So, if $x<x_{n,1}$, $f_n(x)<0$ and if $x>x_{n,n}$, $f_n(x)>0$. In \cite{BracDimRanga} and \cite{Brezinski} we may found the following result about the behaviour of the zeros of  $R_n(x)$.

\begin{theorem}\label{interlace}
\begin{enumerate}
\item The zeros $x_1<\ldots<x_n$ of $R_n$ are real and distinct and at most one of them lies outside $(a,b)$.
\item
\begin{enumerate}
\item If $c>0$, then $x_{n,i}<x_i<x_{n-1,i}$ for $i=1,\ldots,n-1$ and $x_{n,n}<x_n$.
\item If $c<0$, then $x_1<x_{n,1}$ and $x_{n-1,i-1}<x_i<x_{n,i}$ for $i=2,\ldots,n$.
\end{enumerate}
\item If $c<f_n(a)<0$, then $x_1<a$.
\item If $c>f_n(b)>0$, then $b<x_n$.
\item If $f_n(a)<c<f_n(b)$, then $R_n$ has all its zeros in $(a,b)$.
\item Moreover, each $x_i$ is an increasing function of $c$.
\end{enumerate}
\end{theorem}


About the zero location of a real self-reciprocal polynomials of even degree we have the following result, explored by \cite{LAK}.

\begin{theorem} Let $P(z)$ be a real self-reciprocal polynomial of degree $2n$.
\begin{enumerate}
\item $P(z)$ has all its zeros on the unit circle if and only if $C_n(x)$ has all its zeros in $[-1,1]$.
\item If all the zeros of $C_n(x)$ are in $[-1,1]$, written as $x_j=\cos \theta_{j}$, $\theta_{j}\in [0,\pi]$, $j=1,\ldots,n$, then all the zeros $z_{k}$ of $P(z)$, $k=1,\ldots,n$, are given by $z_{k}=e^{\pm i\theta_{k}}$, where
$$
0\leq \theta_{1}<\theta_{2}<\cdots<\theta_{n}\leq \pi.
$$
\end{enumerate}
\end{theorem}
\proof
Firstly, from \eqref{changevariable} follows that
\begin{eqnarray}
z^2+1&=&2xz \label{variable1}\\
z_j+\frac{1}{z_j}&=&2x_j, \mbox{  } j=1,\ldots,n. \label{variable2}
\end{eqnarray}

\begin{enumerate}
\item Suppose that $P(z)$ has all its zeros on the unit circle. As $P(z)$ is a real self-reciprocal polynomial of degree $2n$, its zeros can be arranged in the form $(z_j,\bar{z_j})$, $j=1,\ldots,n$, where $\bar{z_j}=\displaystyle{\frac{1}{z_j}}$ and $|z_j|^2=z_j\bar{z_j}=1$.

We can represent $P(z)$ in the form
\begin{equation}\label{decomposition}
P(z)=p_{2n}\prod_{j=1}^{n}(z-z_j)(z-\bar{z_j})=p_{2n}\prod_{j=1}^{n}(z^2-(z_j+\bar{z_j})z+1).
\end{equation}

Substituting \eqref{variable1} and \eqref{variable2} in \eqref{decomposition} we  have
\begin{equation*}
P(z)=2z^np_{2n}\prod_{j=1}^{n}(x-x_j).
\end{equation*}

Comparing the above equation with \eqref{sumtrans2} follows that
$$C_n(x)=p_{2n}\prod_{j=1}^{n}(x-x_j),$$ 
i.e., $x_j$, $j=1,\ldots,n$,  are the zeros of $C_n(x)$. Furthermore,
$$x_j=\frac{1}{2}\left(z_j+\frac{1}{z_j}\right)=\operatorname{Re} (z_j)$$
and
$$|x_j|=|\operatorname{Re}  (z_j)|\leq |z_j|=1,$$
$j=1,\ldots,n$. Hence $C_n(x)$ has all its zeros in $[-1,1]$.

Now we assume that $C_n(x)$ has all its zeros in $[-1,1]$. From \eqref{sumtrans2}, \eqref{variable1} and \eqref{variable2}, 
\begin{eqnarray*}
P(z)&=&2z^nC_n(z)=2z^np_{2n}\prod_{j=1}^{n}(x-x_j)=p_{2n}\prod_{j=1}^{n}(2xz-2x_jz)\\
&=&p_{2n}\prod_{j=1}^{n}\left(z^2-\left(z_j+\frac{1}{z_j}\right)z+1\right)=p_{2n}\prod_{j=1}^{n}(z-z_j)(z-\bar{z_j}),
\end{eqnarray*}
where $\bar{z_j}=\displaystyle{\frac{1}{z_j}}$, $j=1,\ldots,n$.

Furthermore, from \eqref{variable2},
$$z_j^2-2x_jz_j+1=0 \Leftrightarrow z_j=x_j \pm \left(\sqrt{1-x_j^2}\right)i $$
and
$$|z_j|^2=z_j\bar{z_j}=x_j^2+(1-x_j^2)=1,$$
showing that $|z_j|=1$, $j=1,\ldots,n$.

\item Suppose that all the zeros of $C_n(x)$ are in $[-1,1]$, written as $x_j=\cos \theta_{j}$, $\theta_{j}\in [0,\pi]$, $j=1,\ldots,n$. As  the zeros $z_j$ of $P(z)$ are on the unit circle we can write $z_j=e^{\pm i\gamma_j}$, $0\leq \gamma_j\leq \pi$, $j=1,\ldots,n$. From $x_j=\displaystyle{\frac{1}{2}\left(z_j+\frac{1}{z_j}\right)}$ we have 
$$x_j=\frac{1}{2}(e^{\pm i\gamma_j}+e^{\mp i\gamma_j})=\cos \gamma_j$$
and, consequently, $\theta_j=\gamma_j$, $j=1,\ldots,n$.
\end{enumerate}

\section{Main results}


The first important result of this paper is the following.

\begin{theorem}\label{representacao}
The Chebyshev quasi-orthogonal polynomials of order one generate  the following classes of real self-reciprocal polynomials, given by
\begin{eqnarray*}
P_{\alpha}(z)=\sum_{j=0}^{2n}p_{j,\alpha}z^j ,
\end{eqnarray*}
where $\alpha=T,U,V,W$ (related with each class of Chebyshev polynomials) and
\begin{enumerate}
\item $p_{2n-2,T}=p_{2n-3,T}=\ldots=p_{n,T}=0$ and $p_{2n,T}p_{2n-1,T}\neq 0$;
\item $p_{2n,U}=p_{2n-2,U}=\ldots=p_{n+2,U}=p_{n,U} \neq 0$ and  $p_{2n-1,U}=p_{2n-3,U}=\ldots=p_{n+1,U}\neq 0$;
\item $p_{2n-1,V}=-p_{2n-2,V}=p_{2n-3,V}=\ldots=(-1)^{n}p_{n+1,V}=(-1)^{n-1}p_{n,V}$, $p_{2n,V}\neq 0$ and $p_{2n,V}\neq -p_{2n-1,V}$;
\item $p_{2n-1,W}=p_{2n-2,W}=\ldots=p_{n,W}$,  $p_{2n,W}\neq 0$ and $p_{2n,W}\neq p_{2n-1,W}$;
\end{enumerate}
%
\end{theorem}
\proof
If $p_{2n-2,T}=p_{2n-3,T}=\ldots=p_{n,T}=0$ and $p_{2n,T}p_{2n-1,T}\neq 0$ we have 
$$R_{n,T}(x)=p_{2n,T}T_n(x)+p_{2n-1,T}T_{n-1}(x)=p_{2n,T}\left(T_n(x)+\frac{p_{2n-1,T}}{p_{2n,T}}T_{n-1}(x)\right),$$ 
which  is a first kind Chebyshev quasi-orthogonal polynomial of order one related to $P_T(z)$. In this case, the elements of Theorem \ref{interlace} are $c=c_T=-\displaystyle{\frac{p_{2n-1,T}}{p_{2n,T}}}$, $f_n(1)=f_{n,T}(1)=1$ and $f_n(-1)=f_{n,T}(-1)=-1$. So, 
$$P_T(z)=p_{2n,T}z^{2n}+p_{2n-1,T}z^{2n-1}+p_{2n-1,T}z+p_{2n,T}.$$

In order to obtain a representation of $P_U(z)$, $P_V(z)$ and $P_W(z)$ as a second, third and fourth kinds Chebyshev quasi-orthogonal polynomials of order one, respectively, we have
\begin{eqnarray}
R_{n,U}(x)&=&c_{0,U}U_n(x)+c_{1,U}U_{n-1}(x), \label{quasiU}\\
R_{n,V}(x)&=&c_{0,V}V_n(x)+c_{1,V}V_{n-1}(x),\label{quasiV}\\
R_{n,W}(x)&=&c_{0,W}W_n(x)+c_{1,W}W_{n-1}(x)\label{quasiW},
\end{eqnarray}
where $c_{0,\alpha}c_{1,\alpha}\neq 0$ for $\alpha=U,V,W$.

From \eqref{relationU} and basic manipulations we can rewrite the equation \eqref{quasiU} as
\begin{eqnarray*}
R_{n,U}(x)&=&c_{0,U}\left(2T_n(x)+2T_{n-2}(x)+\ldots+2T_2(x)+T_0(x)\right)\\
&&+c_{1,U}\left(2T_{n-1}(x)+2T_{n-3}(x)+\ldots+2T_1(x)\right).
\end{eqnarray*}
Comparing the coefficients of (\ref{sumtrans2})
with $P_{U}(z)=2z^nR_{n,U}(x)$
follows that 
$$2c_{0,U}=p_{2n,U}=p_{2n-2,U}=\ldots=p_{n+2,U}=p_{n,U}$$
and 
$$2c_{1,U}=p_{2n-1,U}=p_{2n-3,U}=\ldots=p_{n+1,U},$$
obtaining 
$$P_U(z)=p_{2n,U}z^{2n}+p_{2n-1,U}z^{2n-1}+p_{2n,U}z^{2n-2}+\ldots+p_{2n-1,U}z+p_{2n,U}.$$

In this case, $c=c_U=-\displaystyle{\frac{p_{2n-1,U}}{p_{2n,U}}}$, $f_n(1)=f_{n,U}(1)=\displaystyle{\frac{n+1}{n}}$ and $f_n(-1)=f_{n,U}(-1)=-\displaystyle{\frac{n+1}{n}}$.

Following the same idea and using the equations \eqref{relationV} and \eqref{relationW} in \eqref{quasiV} and \eqref{quasiW}, respectively, we obtain  
\begin{eqnarray*}
P_V(z)&=&p_{2n,V}z^{2n}+p_{2n-1,V}z^{2n-1}-p_{2n-1,V}z^{2n-2}+p_{2n-1,V}z^{2n-3}+\ldots\\
&&-p_{2n-1,V}z^2+p_{2n-1,V}z+p_{2n,V}
\end{eqnarray*}
and
$$P_W(z)=p_{2n,W}z^{2n}+p_{2n-1,W}(z^{2n-1}+z^{2n-2}+\ldots+z)+p_{2n,W}.$$

Indeed, $c=c_V=-\left(\displaystyle{\frac{p_{2n-1,V}}{p_{2n,V}}+1}\right)$, $c=c_W=\displaystyle{1-\frac{p_{2n-1,W}}{p_{2n,W}}}$, $f_{n,V}(1)=1$, $f_{n,V}(-1)=-\displaystyle{\frac{2n+1}{2n-1}}$, $f_{n,W}(1)=\displaystyle{\frac{2n+1}{2n-1}}$ and $f_{n,W}(-1)=-1$.

\begin{remark}\label{remark0}
If $p_{2n,V}= -p_{2n-1,V}$ in item 3 or $p_{2n,W}= p_{2n-1,W}$ in item 4  of Theorem \ref{representacao} we have polynomials which coefficients satisfy the conditions of item 2 of Theorem \ref{representacao}. 
\end{remark}

\begin{remark}\label{remark1}
We can observe that when the degree of a real self-reciprocal polynomial $S(z)$ is odd, we have
\begin{eqnarray*}S(z)=s_{2n+1}z^{2n+1}+\ldots+s_{1}z+s_{0}=s_{2n+1}(z+1)P(z),
\end{eqnarray*}
where $P(z)$ is a real self-reciprocal polynomial of degree  $2n$.

Considering 
$$S_{\alpha}(z)=\displaystyle{\sum_{j=0}^{2n+1}s_{j,\alpha}z^j}=s_{2n+1,\alpha}(z+1)P_{\alpha}(z),$$ 
$\alpha = T, U, V, W$, from Theorem \ref{representacao} we find the expressions for the coefficients of $S_{\alpha}(z)$, given by
\begin{enumerate}
\item $s_{2n-2,T}=s_{2n-3,T}=\ldots=s_{n+1,T}=0$, $s_{2n+1,T}s_{2n-1,T}\neq 0$ and $s_{2n}=s_{2n+1}+s_{2n-1}$;
\item $s_{2n,U}=s_{2n-1,U}=\ldots=s_{n+1,U}$ and  $s_{2n+1,U}\neq 0$;
\item $s_{2n-1,V}=s_{2n-2,V}=s_{2n-3,V}=\ldots=s_{n+1,V}=0$ and $s_{2n+1,V}s_{2n,V}\neq 0$;
\item $s_{2n-1,W}=s_{2n-2,W}=\ldots=s_{n+1,W}$, $s_{2n+1,W}s_{2n-1,W}\neq 0$ and \linebreak $\displaystyle{s_{2n,W}=s_{2n+1,W}+\frac{s_{2n-1,W}}{2}}$,
\end{enumerate}
i.e.,
\begin{eqnarray*}
S_T(z)&=&s_{2n+1,T}z^{2n+1}+(s_{2n+1,T}+s_{2n-1,T})z^{2n}+s_{2n-1,T}z^{2n-1}+s_{2n-1,T}z^2\\
&&+(s_{2n+1,T}+a_{2n-1,T})z+s_{2n+1,T},\\
S_U(z)&=&s_{2n+1,U}z^{2n+1}+s_{2n,U}\sum_{j=1}^{2n}z^{j}+s_{2n+1,U},\\
S_V(z)&=&s_{2n+1,V}z^{2n+1}+s_{2n,V}z^{2n}+s_{2n,V}z+s_{2n+1,V},\\
S_W(z)&=&s_{2n+1,W}z^{2n+1}+\left(s_{2n+1,W}+\frac{s_{2n-1,W}}{2}\right)z^{2n}+s_{2n-1,W}\sum_{j=2}^{2n-1}z^{j}\\
&&+\left(s_{2n+1,W}+\frac{s_{2n-1,W}}{2}\right)z+s_{2n+1,W}.
\end{eqnarray*}
\end{remark}

Now, considering  $c_T=-\displaystyle{\frac{p_{2n-1,T}}{p_{2n,T}}}$, $c_U=-\displaystyle{\frac{p_{2n-1,U}}{p_{2n,U}}}$, $c_V=-\left(\displaystyle{\frac{p_{2n-1,V}}{p_{2n,V}}+1}\right)$, $c_W=\displaystyle{1-\frac{p_{2n-1,W}}{p_{2n,W}}}$, $f_{n,T}(1)=f_{n,V}(1)=1$, $f_{n,T}(-1)=f_{n,W}(-1)=-1$,  $f_{n,U}(1)=\displaystyle{\frac{n+1}{n}}$, $f_{n,U}(-1)=-\displaystyle{\frac{n+1}{n}}$,  $f_{n,V}(-1)=-\displaystyle{\frac{2n+1}{2n-1}}$ and  $f_{n,W}(1)=\displaystyle{\frac{2n+1}{2n-1}}$, we can rewrite Theorem \ref{interlace} in the following way.

\begin{lemma}\label{lemma1}
About the zeros of $R_{n,\alpha}(x)$ ($\alpha=T, U, V, W$), we have the following:
\begin{enumerate}
\item The zeros $x_{1,\alpha}<\ldots<x_{n,\alpha}$ of $R_{n,\alpha}$ are real and distinct and at most one of them lies outside $(-1,1)$.
\item
\begin{enumerate}
\item If $c_{\alpha}>0$, then $x_{n,i}<x_{i,\alpha}<x_{n-1,i}$ for $i=1,\ldots,n-1$ and $x_{n,n}<x_{n,\alpha}$.
\item If $c_{\alpha}<0$, then $x_{1,\alpha}<x_{n,1}$ and $x_{n-1,i-1}<x_{i,\alpha}<x_{n,i}$ for $i=2,\ldots,n$.
\end{enumerate}
\item If $c_{\alpha}<f_{n,\alpha}(-1)<0$, then $x_{1,\alpha}<-1$.
\item If $c_{\alpha}>f_{n,\alpha}(1)>0$, then $1<x_{n,\alpha}$.
\item If $f_{n,\alpha}(-1)<c_{\alpha}<f_{n,\alpha}(1)$, then $R_{n,\alpha}$ has all its zeros in $(-1,1)$.
\item Moreover, each $x_{i,\alpha}$ is an increasing function of $c_{\alpha}$.
\end{enumerate}
\end{lemma}

%

The second main result of this paper is the following:

\begin{theorem}\label{main2}About the zeros of $P_{\alpha}(z)$ and  $ S_{\alpha}(z)$ ($\alpha=T, U, V, W$), we have:
\begin{enumerate}
\item Their zeros  are distinct (except in the case that $c_{\alpha}=f_{n,\alpha}(1)$ or $c_{\alpha}=f_{n,\alpha}(-1)$) and at most two of them lie outside the unit circle.
\item If $c_{\alpha}=f_{n,\alpha}(\pm 1)$,  $z=\pm 1$ is zero of multiplicity two of $P_{\alpha}(z)$. For $c_{\alpha}=f_{n,\alpha}(1)$, $z=1$  is zero of multiplicity two of $S_{\alpha}(z)$ and if $c_{\alpha}=f_{n,\alpha}(-1)$, $z=-1$  is zero of multiplicity three of $S_{\alpha}(z)$.
\item If $c_{\alpha}<f_{n,\alpha}(-1)<0$, $P_{\alpha}(z)$ (or $S_{\alpha}(z)$) has two negative zeros and $n-2$ (or $n-1$) distinct zeros on the unit circle. We can represent the negative zeros by $z_{k,\alpha} \in (-\infty,-1)$ and $1/z_{k,\alpha} \in (-1,0)$. 
\item If $c_{\alpha}>f_{n,\alpha}(1)>0$, $P_{\alpha}(z)$ (or $S_{\alpha}(z)$) has two positive zeros and $n-2$ (or $n-1$) distinct zeros on the unit circle. The positive zeros can be represented by $z_{k,\alpha} \in (1,\infty)$ and $1/z_{k,\alpha} \in (0,1)$. 
\item If $f_{n,\alpha}(-1)<c_{\alpha}<f_{n,\alpha}(1)$, $P_{\alpha}(z)$ and $S_{\alpha}(z)$ have all its zeros on the unit circle. They are represented by $z_{k,\alpha}=e^{\pm i\theta_{k,\alpha}}$, with $\theta_{k,\alpha}=\arccos x_{k,\alpha}$, where $x_{k,\alpha}$, $k=1,\ldots,n$, are the zeros of $R_{n,\alpha}(x)$ (for $S_{\alpha}(z)$, we have $x_{n+1,\alpha}=-1$ and $\theta_{n+1,\alpha}=\arccos(-1)=\pi$). Furthermore,
$$
0\leq \theta_{1,\alpha}<\theta_{2,\alpha}<\cdots<\theta_{n,\alpha}\leq \pi,
$$
$\theta_{n+1,\alpha}=\pi$ (odd case) and $\theta_{j,\alpha}$, $j=1,\ldots,n$  are increasing functions of $c_{\alpha}$.
\end{enumerate}
\end{theorem}
\proof
The proof follows direct from (\ref{changevariable}), Theorem \ref{representacao}, Remark \ref{remark1} and Lemma \ref{lemma1}. 

About the multiplicities,
$$P_{\alpha}(z)=0  \Leftrightarrow R_{n,\alpha}(x)=0 \Leftrightarrow c_{\alpha}=f_{n,\alpha}(x).$$

Furthermore, 
$$P'_{\alpha}(z)=0  \Leftrightarrow z=\pm 1 \mbox{ (consequently, } x=\pm 1)$$
and $P''_{\alpha}(\pm 1)\neq 0.$

%
%


%



\end{document}